\documentclass[twocolumn]{autart}    

\usepackage{graphicx}          
\usepackage{color}
\usepackage[dvips]{epsfig}
\usepackage{epstopdf}

\begin{document}

\begin{frontmatter}

\title{Generalized Lyapunov criteria on finite-time stability of stochastic nonlinear systems\thanksref{footnoteinfo}} 

\thanks[footnoteinfo]{This paper was not presented at any IFAC
meeting, and was supported in part by the Natural Science Foundation of China under Grant 61573006. Corresponding author: J. Yin}

\author[Yu]{Xin Yu}\ead{xyu@ujs.edu.cn},    
\author[Yin]{Juliang Yin}\ead{yin$\_$juliang@hotmail.com},               
\author[Khoo]{Suiyang Khoo}\ead{sui.khoo@deakin.edu.au}  

\address[Yu]{School of Electrical and Information Engineering, Jiangsu University, Zhenjiang, China}  
\address[Yin]{School of Economics and Statistics, Guangzhou University, Guangzhou, China}             
\address[Khoo]{School of Engineering, Deakin University, VIC, Australia}        

\begin{keyword}                           
Finite-time stability; Generalized Lyapunov theorem; Multiple Lyapunov functions; Stochastic nonlinear systems.               
\end{keyword}                             

\begin{abstract}                          
This paper considers the problem of finite-time stability for stochastic nonlinear systems. A new Lyapunov theorem of stochastic finite-time stability is proposed, and an important corollary is obtained. Some comparisons with the existing results are given, and it shows that this new Lyapunov theorem not only is a generalization of classical stochastic finite-time theorem, but also reveals the important role of white-noise in finite-time stabilizing stochastic systems. In addition, multiple Lyapunov functions-based criteria on stochastic finite-time stability are presented, which further relax the constraint of the infinitesimal generator $\mathcal{L}V$.  Some examples are constructed to show significant features of
the proposed theorems. Finally, simulation results are presented to
demonstrate the theoretical analysis.
\end{abstract}

\end{frontmatter}

\section{Introduction}
Stochastic stability has been one of the most fundamental research topics of controlled systems modeled by stochastic differential equation in the past few decades, and we here mention \cite{Kusher}, \cite{Arnold72}, \cite{Khasminskii}, \cite{Krstic2}, \cite{deng01}, \cite{mao94}, \cite{mao07}, \cite{ito15} and \cite{zhao16} among others.

In classical stochastic stability theory, asymptotic stability in probability, $p$-order moment asymptotic stability, and
almost sure asymptotic stability are often considered. These three types
of stability describe the asymptotic
behavior of the trajectories of a stochastic system as time goes
to infinity.  In many applications, however, it is desirable
that a stochastic system possesses the property that its trajectories
converge to a Lyapunov stable equilibrium state
in finite time rather than merely asymptotically.

To develop the theory of finite-time stability of deterministic systems \cite{bhat00} to the stochastic case, the notions and Lyapunov criteria of stochastic finite-time stability were introduced separately in \cite{yin11} and \cite{chen10}. Based on the presented stochastic finite-time stability theory, finite-time stabilization controllers of stochastic nonlinear systems by state or output feedback are designed in \cite{khoo13}, \cite{yin15}, \cite{wang15} and \cite{zha14} for example. Properties of finite-time stable stochastic systems are further discussed in \cite{yin2015}. Recently, finite-time stability of homogeneous stochastic nonlinear systems is studied in \cite{yin17}.

In those existing papers, the stochastic finite-time stability is required that the infinitesimal generator $\mathcal{L}V$ satisfies $\mathcal{L}V\leq-cV^{\gamma}$ with $0<\gamma<1$ and $c>0$. So far, to the best of our
knowledge, there are not any papers that could demonstrate whether the stochastic finite-time
stability still holds or not, if this condition is not fulfilled.

In this paper, our target
is to establish new Lyapunov criteria of stochastic finite-time stability under more general conditions. The main contributions of this paper are as follows: A new Lyapunov theorem on finite-time stability of stochastic nonlinear systems is proved, and an important corollary follows directly; By comparing the new Lyapunov theorem with the previous results of stochastic finite-time stability, it is shown that this generalized criterion not only relaxes the constraint on infinitesimal generator $\mathcal{L}V$,  but also reveals the important applications of white-noise in finite-time stabilizing a system; In addition, multiple Lyapunov functions-based criteria of stochastic finite-time stability are presented, which further relax the constraint of the infinitesimal generator $\mathcal{L}V$; Some examples are constructed to show the significant features of
our results.

The rest of the paper is organized as follows. The mathematic preliminaries are given in Section 2. A new Lyapunov theorem of stochastic finite-time stability is presented in Section 3. In Section 4, we derive an important corollary and discuss the comparisons with the existing results.  In Section 5, multiple Lyapunov functions-based criteria of stochastic finite-time stability are presented. Section 6 gives some simulation
results to illustrate the theoretical results. Finally, concluding
remarks are given in Section 7.

\textbf{Notations:} \ $\mathrm{R}_+$ stands for the set of all nonnegative real
numbers, $\mathrm{R}^n$ is the $n$-dimensional Euclidean space, $\mathrm{R}^{n\times
m}$ is the space of real $n\times m$-matrices. $|x|$ is the usual
Euclidean norm of a vector $x$. $A^T$ denotes the transpose of matrix
$A$.
$\mathrm{Tr}\{A\}$ is its
trace when $A$ is a square matrix.  $C^2$ denotes the family of all functions with continuous
second partial derivatives. A random variable $\xi\in L^1$ means that $\mathrm{E}|\xi|<\infty$. $\mu:\mathrm{R}_+\rightarrow \mathrm{R}_+$ is a $\mathcal{K}_\infty$ class function means that
it is continuous, strictly increasing,
$\mu(0)=0$ and $\lim_{s\rightarrow\infty}\mu(s)=\infty$.
 $a\wedge b=\min\{a,b\}$.

\section{Preliminary results}

In this paper, we consider a stochastic nonlinear system modeled by the following stochastic differential
equation:
\begin{eqnarray} \label{a1}
dx=f(x)dt+g(x)dB(t), \quad x(0)=x_0\in \mathrm{R}^n,
\end{eqnarray}
where $x\in \mathrm{R}^n$ is the system state; $B(\cdot)$ is an $m$-dimensional standard Brownian motion defined on a complete
probability space $(\Omega, \mathcal{F},
 \mathrm{P})$;
 $f:\mathrm{R}^n\rightarrow \mathrm{R}^n$  and $g: \mathrm{R}^n\rightarrow \mathrm{R}^{n\times m}$ are continuous in $x$ and satisfying
$f(0)=0$ and $g(0)=0$, which implies that (\ref{a1}) has a trivial zero solution.

As discussed in \cite{yin17}, in general, we are interested in having a unique solution in forward time
for a stochastic differential equation. However, it is generally difficult
to ensure such a property for a stochastic differential equation
without locally Lipschitz continuous coefficients. Actually, it has
been pointed out in \cite{yin11}, a finite-time stable stochastic nonlinear
system does not have locally Lipschitz continuous coefficients. Hence, it suffices to require the existence
of solutions for a stochastic nonlinear system either in the strong sense or in the weak sense when studying
finite-time stochastic stability.

\textbf{Remark 1: } A weak solution to system (\ref{a1}) can be well-defined, and its precise definition can be found in (\cite{rogers94}, p.149). In fact, a strong solution is of course also a weak solution. It is also pointed out in \cite{rogers94} that the concept of weak solutions is appropriate for control problems.

Let $\mathcal {L}V$ denote infinitesimal generator of a $C^2$ function
$V:\mathrm{R}^n\rightarrow\mathrm{R}$ along stochastic differential equation (\ref{a1})
with the definition of
\begin{eqnarray}\label{a2}
\mathcal {L}V(x)=\frac{\partial V(x)}{\partial x}f(x)+\frac{1}{2}\mbox{Tr}\left\{g^T(x)\frac{\partial^2V(x)}{\partial
x^2}g(x)\right\},
\end{eqnarray}
where $\frac{\partial V}{\partial x}$ denotes the gradient of $V$ (written as a row vector), and $\frac{\partial^2V}{\partial
x^2}$ denotes the Hessian of $V$.

The following lemma (see Lemma 2.1 of \cite{yin15}) gives an existence
result of a weak solution to system (\ref{a1}).

\textbf{Lemma 1 \cite{yin15}:} Suppose that there exists a nonnegative $C^2$ function
$V:\mathrm{R}^n\rightarrow\mathrm{R}_+$, which is radially unbounded, that is,
$\lim_{|x|\rightarrow\infty}V(x)=\infty$. If the infinitesimal generator of $V$
with respect to (\ref{a1}) satisfies $\mathcal {L}V(x)\leq0$, then (\ref{a1}) has a regular continuous solution (in the weak sense) for any initial data.

A regular solution means that the solution has no finite explosion time with probability one. The detailed definition of regular solution can be founded in \cite{Khasminskii}. The next lemma is one of the well-known Doob's Optional-Sampling Theorem for continuous nonnegative supermartingales, which can be found in (\cite{rogers94'}, p.189, (77.5)) and is useful in later analysis.

\textbf{Lemma 2 \cite{rogers94'}:} Suppose that $X(t)$ is a continuous nonnegative supermartingale with respect to a filtration $\{\mathcal{F}_t\}_{t\geq0}$. Let $S$ and $T$ be stopping times with $S\leq T$. Then $X(T)\in L^1$, and
\begin{eqnarray}\label{a3}
\mathrm{E}(X(T)|\mathcal{F}_S)\leq X(S).
\end{eqnarray}

\section{Generalized stochastic finite-time stability theorem}
In this section, we first review and refine the definition of stochastic finite-time
stability introduced in \cite{yin11,yin15}. Then a new Lyapunov theorem on finite-time
stability of stochastic nonlinear systems will be given, and an important corollary is derived as well.

\textbf{Definition 1:} The trivial zero solution of
(\ref{a1}) is said to be stochastically finite-time stable, if the stochastic
system admits a solution (either in the strong sense or in the weak sense) for any initial data $x_0\in \mathrm{R}^n$,  and the following properties hold:\\
(i) Finite-time attractiveness in probability: For every initial value $x_0\in \mathrm{R}^n\setminus\{0\}$, and any solution $x(t;x_0)$, the first hitting
time of $x(t;x_0)$, i.e., $\tau_{x_0}=\inf\{t\geq0: x(t;x_0)=0\}$, called stochastic settling time, is
finite almost surely, that is $\mathrm{P}(\tau_{x_0}<\infty)=1$; moreover,
\begin{eqnarray}\label{a4}
x(t+\tau_{x_0};x_0)=0, \  \forall t\geq0, \ \hbox{a.s.}
\end{eqnarray}
(ii) Stability in probability: For any solution $x(t;x_0)$, and every pair of $\varepsilon\in(0,1)$ and $r>0$,
there exists a $\delta(\varepsilon,r)>0$ such that
\begin{eqnarray}\label{a5}
\mathrm{P}(|x(t;x_0)|< r \ \hbox{for all} \ t>0)\geq1-\varepsilon
\end{eqnarray}
whenever $|x_0|<\delta$.

\textbf{Remark 2:} The finite-time attractiveness in probability defined here states
that any trajectories of a stochastic system will not only reach the origin
in finite time, but also stay at the origin for ever after the stochastic settling time almost surely. So the origin is both an equilibrium point and an absorbing state.

\textbf{Remark 3:}  The stability in probability is equivalent to that: For any solution $x(t;x_0)$, and any
$r>0$,
\begin{eqnarray}\label{a6}
\lim_{x_0\rightarrow0}\mathrm{P}\left(\sup_{t\geq0}|x(t;x_0)|\geq r\right)=0
\end{eqnarray}
holds, which will be used in the following analysis.

Now, it is ready to state a new Lyapunov theorem on stochastic finite-time stability.

\textbf{Theorem 1:} For system (\ref{a1}), if there exists a $C^2$ positive definite and radially unbounded function $V:\mathrm{R}^n\rightarrow\mathrm{R}_+$, a positive constant $c>0$, such that
\begin{eqnarray}
 \mathcal {L}V(x)\leq0,  \  \forall x\in \mathrm{R}^n, \label{a7}
\end{eqnarray}
and
\begin{eqnarray}\label{a8}
&&K(V(x))\left[cK(V(x))+\mathcal {L}V(x)\right] \nonumber\\
&&\leq \frac{K{'}(V(x))}{2}\left|\frac{\partial V}{\partial x}g(x)\right|^2, \  \forall x\in \mathrm{R}^n\setminus\{0\},
\end{eqnarray}
where $K:\mathrm{R}_+\rightarrow \mathrm{R}_+$ is a continuous differentiable function with the derivative $K'(s)\geq0$ and $K(s)>0$ for any $s>0$ and
\begin{eqnarray}\label{a9}
\int_{0}^\epsilon \frac{ds}{K(s)}<\infty, \  \forall\epsilon>0,
\end{eqnarray}
then the trivial solution of (\ref{a1}) is stochastically finite-time stable, and stochastic settling time satisfies
\begin{eqnarray}\label{a9.1}
\mathrm{E}\tau_{x_0}\leq \frac{1}{c}\int_{0}^{V(x_0)}\frac{ds}{K(s)}.
\end{eqnarray}


\textbf{Remark 4:} It is easy to know the following  functions have the same properties as $K$ in Theorem 1:
\begin{eqnarray*}
K_1(s)&=&s^{\gamma}, \quad 0<\gamma<1, s\ge 0,\\
K_2(s)&=&s^{\gamma}+s^{\alpha}, \quad  0<\gamma<1, \alpha\geq1, s\ge 0.
\end{eqnarray*}

Before giving the proof of Theorem 1, we need the following lemma.

\textbf{Lemma 3:} Let $(X(t), \mathcal{F}_t; 0\leq t<\infty)$ be a continuous nonnegative supermartingale and $\tau=\inf\{t\geq0; X(t)=0\}$. If $\mathrm{P}(\tau<\infty)=1$, then
\begin{eqnarray}\label{a11}
X(t+\tau)=0, \ \forall t\geq0, \ \hbox{a.s.}
\end{eqnarray}

\textbf{Proof:} Since $X(t)$ is continuous,
$\tau=\inf\{t\geq0; X(t)=0\}$ is a stopping time. So, for any real constant $t\geq0$, $t+\tau$ is also a stopping time. By Lemma 2, we have $X(t+\tau)\in L^1$, and
\begin{eqnarray}\label{a12}
\mathrm{E}(X(t+\tau)|\mathcal{F}_{\tau})\leq X(\tau), \ \forall t\geq0.
\end{eqnarray}
Taking expectation on both sides of (\ref{a12}), one gets
\begin{eqnarray}\label{a13}
\mathrm{E}X(t+\tau)\leq \mathrm{E}X(\tau)=0, \ \forall t\geq0,
\end{eqnarray}
which together with the nonnegativity of $X(t)$ leads to
\begin{eqnarray}\label{a14}
\mathrm{E}X(t+\tau)=0, \ \forall t\geq0.
\end{eqnarray}
Therefore, this implies that (\ref{a11}) holds, which completes the proof.

Now, we can give the detailed proof of Theorem 1.

\textbf{Proof of Theorem 1:} From (\ref{a7}) and Lemma 1, it leads to that for
each $x_0\in \mathrm{R}^n$, there exists a regular continuous solution $x(t;x_0)$
to (\ref{a1}). Since $\mathcal{L}V(x)\leq0$ and $V(x)\geq0$, $V_t=V(x(t;x_0))$ is a nonnegative continuous supermartingale with augmented filtration $\{\mathcal{F}_t\}_{t\geq0}$ satisfying
the usual conditions.

Since $V$ is positive definite and radially unbounded, from \cite{Khalil}, there exists a $\mathcal{K}_\infty$ class function $\mu$ such that
\begin{eqnarray}\label{a15'}
\mu(|x|)\leq V(x).
\end{eqnarray}
By a supermartingale inequality (\cite{mao07}, p.13, Theorem 3.6), for any $r>0$ and any natural number $n$, one has
\begin{eqnarray}\label{a15}
\mathrm{P}\left(\sup_{0\leq t\leq n}V_t\geq\mu(r)\right)\leq \frac{\mathrm{E}V_0+\mathrm{E}V_n}{\mu(r)}\leq\frac{2V(x_0)}{\mu(r)}.
\end{eqnarray}
From $\mu(|x|)\leq V(x)$, it leads to
\begin{eqnarray}\label{a16}
\left\{\sup_{0\leq t\leq n}\mu(|x(t;x_0)|)\geq\mu(r)\right\}\subseteq \left\{\sup_{0\leq t\leq n}V_t\geq\mu(r)\right\},
\end{eqnarray}
which with (\ref{a15}) and $\mu$ being $\mathcal{K}_\infty$ function implies that
\begin{eqnarray}\label{a17}
&&\mathrm{P}\left(\sup_{0\leq t\leq n}|x(t;x_0)|\geq r\right) \nonumber\\
&&=\mathrm{P}\left\{\sup_{0\leq t\leq n}\mu(|x(t;x_0)|)\geq\mu(r)\right\} \nonumber\\
&&\leq \frac{2V(x_0)}{\mu(r)}.
\end{eqnarray}
By monotonic convergence theorem, we have
\begin{eqnarray}\label{a18}
&&\mathrm{P}\left(\sup_{t\geq 0}|x(t;x_0)|\geq r\right) \nonumber\\
&&=\lim_{n\rightarrow\infty}\mathrm{P}\left(\sup_{0\leq t\leq n}|x(t;x_0)|\geq r\right) \nonumber\\
&&\leq \frac{2V(x_0)}{\mu(r)}.
\end{eqnarray}
So, for any $r>0$, by the continuity of $V$, one has
\begin{eqnarray}\label{a19}
\lim_{x_0\rightarrow0}\mathrm{P}\left(\sup_{t\geq0}|x(t;x_0)|\geq r\right)=0.
\end{eqnarray}
By Remark 3, the stability in probability follows.

We now turn our attention to proving the finite-time attractiveness in probability. Define a positive definite function
\begin{eqnarray}\label{a20}
F(x)=\int_{0}^{V(x)}\frac{ds}{K(s)},
\end{eqnarray}
which can be verified that it is $C^2$ in $\mathrm{R}^n\setminus \{0\}$. For any initial value $x_0\in \mathrm{R}^n\setminus\{0\}$, we define stopping times
\begin{eqnarray}\label{a21}
\tau_{k}=\inf\left\{t\geq 0; |x(t;x_0)|\notin \left(\frac{1}{k},k\right)\right\},
\end{eqnarray}
\begin{eqnarray}\label{a21'}
\tau_{1k}=\inf\left\{t\geq 0; |x(t;x_0)|\in \left[0,\frac{1}{k}\right]\right\},
\end{eqnarray}
and
\begin{eqnarray}\label{a21''}
\tau_{2k}=\inf\left\{t\geq 0; |x(t;x_0)|\in \left[k,\infty\right)\right\},
\end{eqnarray}
with $\inf{\emptyset}=\infty$ and nature numbers $k\in \{2,3,4,\cdots\}$. It is clear
that $\tau_k$, $\tau_{1k}$ and $\tau_{2k}$ are all increasing stopping time sequences, and $\tau_{k}=\tau_{1k}\wedge\tau_{2k}$. We now set $\tau_\infty=\lim_{k\rightarrow\infty}\tau_k$, and $\tau_{i\infty}=\lim_{k\rightarrow\infty}\tau_{ik}$, for $i=1,2$. Since the solution $x(t;x_0)$ is regular, $\tau_{2\infty}=\infty$ a.s., and therefor $\tau_{\infty}=\tau_{1\infty}$ a.s..

The function $F(x)$ is twice continuously differentiable in the domain $\frac{1}{k}< |x|< k$ for any $k$. Applying It\^{o}'s formula in this domain, we get
\begin{eqnarray}\label{a22}
F(x(t\wedge\tau_{k};x_0))&=&F(x_0)+\int_{0}^{t\wedge\tau_{k}}\mathcal{L}F(x(s;x_0))ds \nonumber\\
&&+\int_{0}^{t\wedge\tau_{k}}\frac{1}{K(V)}\frac{\partial V}{\partial x}g(x(s;x_0))dB(s).
\end{eqnarray}
Noting that $\frac{1}{K(V)}\frac{\partial V}{\partial x}g(x)$ is bounded in the domain $\frac{1}{k}< |x|< k$,  we have
\begin{eqnarray}\label{a23'}
\int_{0}^{t}\mathrm{E}\left(\frac{I{\{s\leq\tau_{k}\}}}{K^2(V)}\left|\frac{\partial V}{\partial x}g(x(s;x_0))\right|^2\right)ds<\infty,
\end{eqnarray}
where $I\{\cdot\}$ denotes the indicator function, which together with
\begin{eqnarray}\label{a23}
&&\int_{0}^{t\wedge\tau_{k}}\frac{1}{K(V)}\frac{\partial V}{\partial x}g(x(s;x_0))dB(s) \nonumber\\
&&=\int_{0}^{t}I{\{s\leq\tau_{k}\}}\frac{1}{K(V)}\frac{\partial V}{\partial x}g(x(s;x_0))dB(s)
\end{eqnarray}
implies that $\mathrm{E}(\int_{0}^{t\wedge\tau_{k}}\frac{1}{K(V)}\frac{\partial V}{\partial x}g(x(s;x_0))dB(s))=0$. Taking expectation on both sides of (\ref{a22}), we have
\begin{eqnarray}\label{a24}
\mathrm{E}F(x(t\wedge\tau_{k};x_0))&=&F(x_0)+\mathrm{E}\int_{0}^{t\wedge\tau_{k}}\mathcal{L}F(x(s;x_0))ds.
\end{eqnarray}
In the domain $\frac{1}{k}< |x|< k$, it can be verified that
\begin{eqnarray}\label{a25}
\mathcal{L}F(x)=\frac{\mathcal{L}V(x)}{K(V(x))}-\frac{K'(V(x))}{2K^2(V(x))}\left|\frac{\partial V}{\partial x}g(x)\right|^2.
\end{eqnarray}
By condition (\ref{a8}), it is obvious that in this domain $\mathcal{L}F(x)\leq-c$, which together with (\ref{a24}) and $F\geq0$ leads to
$\mathrm{E}(t\wedge\tau_{k})\leq \frac{F(x_0)}{c}$. Letting $t,k\rightarrow \infty$, using Fatou lemma and $\tau_\infty=\tau_{1\infty}$ a.s., one has
\begin{eqnarray}\label{a26}
\mathrm{E}\tau_{1\infty}=\mathrm{E}\tau_\infty \leq \frac{F(x_0)}{c}.
\end{eqnarray}
It is clear that the first hitting time $\tau_{x_0}=\inf\{t\geq0: x(t;x_0)=0\}=\lim_{k\rightarrow\infty}\tau_{1k}=\tau_{1\infty}$ a.s. Therefore,
\begin{eqnarray}\label{a27}
\mathrm{E}\tau_{x_0}\leq \frac{F(x_0)}{c}=\frac{1}{c}\int_0^{V(x_0)}\frac{ds}{K(s)}<\infty,
\end{eqnarray}
which implies that $\mathrm{P}(\tau_{x_0}<\infty)=1$.

Recalling that $V_t=V(x(t;x_0))$ is a continuous nonnegative supermartingale, using $V_{\tau_{x_0}}=V(x(\tau_{x_0};x_0))=0$ and Lemma 3, we have
\begin{eqnarray}\label{a28}
V_{t+\tau_{x_0}}=0, \ \forall t\geq0, \ \hbox{a.s.}
\end{eqnarray}
From $\mu(|x|)\leq V(x)$, it follows that $x(t+\tau_{x_0};x_0)=0$ a.s., for any $t\geq0$. Here we complete the proof.

\textbf{Remark 5:} Following the same proof process above, we can see that Theorem 1 still holds for nonautonomous stochastic nonlinear systems \begin{eqnarray} \label{a28'}
dx=f(t,x)dt+g(t,x)dB(t)
\end{eqnarray}
with an additional assumption that system (\ref{a28'}) admits a regular continuous solution (either in the strong sense or in the weak sense) for any initial state $x_0\in \mathrm{R}^n$.

\textbf{Remark 6:}  If system (\ref{a1}) degenerates into a deterministic system, i.e., the diffusion dynamic $g(x)=0$ in (\ref{a1}), the condition (\ref{a8}) turns into
\begin{eqnarray}\label{a10}
\dot{V}(x)\leq -cK(V(x)), \ \hbox{for} \ \forall x\in \mathrm{R}^n\setminus\{0\},
\end{eqnarray}
and Theorem 1 then reduced to the corresponding Lyapunov theorem of deterministic system for finite-time stability \cite{Moulay}.

\textbf{Remark 7:}  If condition (\ref{a9}) is replaced by $\int_{0}^{\infty}\frac{ds}{K(s)}\leq C$, where $C$ is a positive constant, then the stochastic settling time $\tau_{x_0}$ satisfies  $\mathrm{E}\tau_{x_0}\leq C/c$, and its upper bound is independent of the initial state $x_0$.  Such a function can be chosen as $K(s)=s^{\gamma}+s^{\alpha}$ with  $0<\gamma<1$ and $\alpha>1$.

\section{Comparisons with the existing results}
Let us first recall the existing results on the stochastic
finite-time stability \cite{yin11,yin15}, and take the classical result in \cite{yin11} as a theorem.

\textbf{Theorem 2 \cite{yin11}}: For system (\ref{a1}), If there
exists a $C^2$ function $V:\mathrm{R}^n\rightarrow \mathrm{R}_{+}$, $\mathcal{K}_{\infty}$ class functions $\mu_1$ and $\mu_2$,
positive real numbers $c>0$ and $0<\gamma<1$, such that for all $x\in \mathrm{R}^n$,
\begin{eqnarray}
&&\mu_1(|x|)\leq V(x) \leq \mu_2(|x|), \label{a29'} \\
&&\mathcal {L}V(x)\leq-c(V(x))^{\gamma},  \label{a30'}
\end{eqnarray}
 then the trivial solution of (\ref{a1}) is stochastically finite-time stable.

To see the important contributions of this paper, let us first state a useful
corollary that follows from Theorem 1 directly.

\textbf{Corollary 1:} For system (\ref{a1}), if there exists a $C^2$ positive definite and radially unbounded function $V: \mathrm{R}^n\rightarrow\mathrm{R}_+$, a positive constant $c>0$, such that
\begin{eqnarray}
 \mathcal {L}V(x)\leq0,  \ \forall x\in \mathrm{R}^n,  \label{a29}
\end{eqnarray}
and
\begin{eqnarray}\label{a30}
V\cdot\left(cV^{\gamma}+\mathcal {L}V\right) \leq \frac{\gamma}{2}\left|\frac{\partial V}{\partial x}g\right|^2,
 \ \forall x\in \mathrm{R}^n\setminus\{0\},
\end{eqnarray}
where $0<\gamma<1$ and the argument $x$ is omitted here without ambiguity,
then the trivial solution of (\ref{a1}) is stochastically finite-time stable, and stochastic settling time satisfies
\begin{eqnarray}\label{a30.1}
\mathrm{E}\tau_{x_0}\leq \frac{1}{c(1-\gamma)}(V(x_0))^{1-\gamma}.
\end{eqnarray}

\textbf{Proof:} Letting $K(s)=s^{\gamma}$, $0<\gamma<1$, in Theorem 1, we have Corollary 1 directly.

Let us explain the significant features of this corollary from the following two aspects.

\textbf{(I)} In the classical Theorem 2, $\mathcal{L}V$ is required to be not only negative definite, but also not greater than a kind of functions $-cV^{\gamma}$ with $0<\gamma<1$. As far as we know, there is not a paper so far that shows whether the stochastic finite-time
stability holds or not if this condition does not hold, but
our Corollary 1 gives a positive answer. In fact, we see from condition
(\ref{a30}) that $\mathcal{L}V$ may be not negative definite somewhere (see the examples below
for an explicit support) but yet the corollary shows that the system may
still be stochastically finite-time stable.

\textbf{(II)} We see clearly that if (\ref{a30'}) is satisfied, (\ref{a30}) must be satisfied but not conversely.
It is the term $|\frac{\partial V}{\partial x}g|^2$ that makes
condition (\ref{a30}) be satisfied much more easily than condition (\ref{a30'}).
So Corollary 1 has already enabled us to construct Lyapunov functions
more easily in applications. Note furthermore that the term $|\frac{\partial V}{\partial x}g|^2$
is connected with the diffusion coefficient $g(x)$, so our result reveals the
important role of white-noise in finite-time stabilizing a stochastic system.

\textbf{Example 1:} Consider a one-dimensional stochastic nonlinear
system in the form
\begin{eqnarray} \label{a31}
dx=f(x)dt+g(x)dB(t), \quad x_0\neq0,
\end{eqnarray}
where
\begin{eqnarray} \label{a32}
&&f(x)=-c_1x^{p}-c_2x^{\beta_1}; \ c_1\geq0, \ c_2>0, \ 1>\beta_1=\frac{p_1}{q_1}, \nonumber \\
&&g(x)=c_3x^{\beta_2}; \ c_3\neq0, \ 1>\beta_2>\frac{1}{2},
\end{eqnarray}
and $p$, $p_1$ and $q_1$ are positive odd numbers.
Consider a $C^2$
Lyapunov function $V(x)=|x|^{\alpha}$ with $\alpha\geq2$. It is not hard to compute
\begin{eqnarray} \label{a33}
\mathcal{L} V(x)&=&-c_1\alpha|x|^{\alpha+p-1}-c_2\alpha|x|^{\alpha+\beta_1-1} \nonumber\\
&&+\frac{1}{2}c_3^2\alpha(\alpha-1)|x|^{\alpha+2\beta_2-2}.
\end{eqnarray}
Let us analyze this example from three concrete cases.

\textbf{Case 1.} If the parameters satisfy $c_1=0$, $2\beta_2=\beta_1+1$ and $c_2=\frac{1}{2}(\alpha-1)c_3^2$,
we have
\begin{eqnarray} \label{a34}
\mathcal{L} V(x)=\mathcal{L} |x|^{\alpha}\equiv0.
\end{eqnarray}
Meanwhile, it is not hard to verify that the condition (\ref{a30}) is satisfied with
$\gamma=\frac{\alpha+2\beta_2-2}{\alpha}$ and $c\leq \frac{1}{2}\gamma c_3^2\alpha^2$.  Clearly $0<\gamma<1$, and hence, the system (\ref{a31}) in this case is stochastically finite-time stable by Corollary 1 even though $\mathcal{L} V(x)=0$.

\textbf{Case 2.} If the parameters satisfy $c_1>0$, $p=1$, $2\beta_2=\beta_1+1$ and $c_2=\frac{1}{2}(\alpha-1)c_3^2$,
we have
\begin{eqnarray} \label{a35}
\mathcal{L} V(x)=-c_1\alpha V(x)=:-c_0V(x).
\end{eqnarray}
The condition (\ref{a30}) is also satisfied with
$\gamma=\frac{\alpha+2\beta_2-2}{\alpha}$ and $c\leq \frac{1}{2}\gamma c_3^2\alpha^2$. Hence, the system (\ref{a31}) in this case is still stochastically finite-time stable by Corollary 1 even though $\mathcal{L} V(x)=-c_0V(x)$.

\textbf{Case 3.} If the parameters satisfy $c_1>0$, $p=3$, $2\beta_2=\beta_1+1$, $c_2=\frac{1}{2}(\alpha-1)c_3^2$ and $\alpha=2$,
we have
\begin{eqnarray} \label{a35.1}
\mathcal{L} V(x)=\mathcal{L} |x|^2=-2c_1(V(x))^2.
\end{eqnarray}
The condition (\ref{a30}) is also satisfied with
$\gamma=\frac{\alpha+2\beta_2-2}{\alpha}$ and $c\leq \frac{1}{2}\gamma c_3^2\alpha^2$. Hence, the system (\ref{a31}) in this case is still stochastically finite-time stable by Corollary 1 even though $\mathcal{L} V(x)=-2c_1 (V(x))^2$.

\textbf{Example 2:} Let us consider a 2-dimensional stochastic nonlinear
system of the form
\begin{eqnarray} \label{a36}
&&dx_1=f_1(x_1,x_2)dt+g_1(x_1,x_2)dB_1(t), \nonumber\\
&&dx_2=f_2(x_1,x_2)dt+g_2(x_1,x_2)dB_2(t),
\end{eqnarray}
where $B_1(t)$ and $B_2(t)$ are two mutually independent Brownian
motions, and the coefficients are expressed as
\begin{eqnarray} \label{a37}
&&f_1(x)=-\frac{1}{8}x_1^{\frac{1}{3}}+x_2, \quad g_1(x)=\frac{1}{2}\hbox{sign}(x_1)|x_1|^{\frac{2}{3}},\nonumber\\
&&f_2(x)=-x_1-\frac{1}{8}x_2^{\frac{1}{3}}, \quad g_2(x)=\frac{1}{2}\hbox{sign}(x_2)|x_2|^{\frac{2}{3}}.
\end{eqnarray}
By using a Lyapunov function $V(x)=x^2_1+ x^2_2$, it is easy to verify that
\begin{eqnarray} \label{a38}
\mathcal{L} V(x)=2x_1\cdot f_1+2x_2\cdot f_2+g_1^2+g_2^2\equiv0.
\end{eqnarray}
However, we can verify condition (\ref{a30}) is satisfied with $\gamma=\frac{2}{3}$ and $c\leq\frac{\gamma}{2}$. In fact, using an elementary
inequality in (\cite{huang}, Lemma 2.3), that is,
\begin{eqnarray} \label{a39}
(a^2+b^2)^{\frac{1+\gamma}{2}}\leq |a|^{1+\gamma}+|b|^{1+\gamma}, \quad 0<\gamma<1,
\end{eqnarray}
we have condition (\ref{a30}) holds, i.e., the following inequality holds
\begin{eqnarray} \label{a40}
c(x_1^2+x_2^2)^{1+\gamma} \leq \frac{\gamma}{2}
(|x_1|^{1+2/3}+|x_2|^{1+2/3})^2,
\end{eqnarray}
 with $\gamma=\frac{2}{3}$ and $c\leq\frac{\gamma}{2}$. So by Corollary 1, system (\ref{a36}) is stochastically finite-time stable even though $\mathcal{L} V(x)=0$.

 \textbf{Remark 8:} Examples 1 and 2 give an explicit illustration about the significant features of Corollary 1 or Theorem 1, and show that stochastic nonlinear system may still be finite-time stable even though  $\mathcal{L}V=0$ or $\mathcal{L}V=-c_0V^{\gamma}$ with $c_0>0$ and $\gamma\geq1$.  In addition, from the Case 2 of Example 1,  we can see that Eq.(4.3) in \cite{yin11} seems to be a necessary condition on stochastic finite-time instability theorem.

\section{Multiple Lyapunov functions-based criteria of stochastic finite-time stability}
In this section, we shall develop Theorem 1 by the use of multiple Lyapunov functions, and obtain the following stochastic finite-time stability criteria.

\textbf{Theorem 3:} For system (\ref{a1}), suppose that there exists a $C^2$ positive definite and radially unbounded function $U:\mathrm{R}^n\rightarrow\mathrm{R}_+$ such that
\begin{eqnarray}\label{a41}
 \mathcal {L}U(x)\leq0,  \  \forall x\in \mathrm{R}^n.
\end{eqnarray}
Furthermore, if there exists a $C^2$ positive definite function $V:\mathrm{R}^n\rightarrow\mathrm{R}_+$, a positive constant $c>0$ and
a continuous differentiable function $K:\mathrm{R}_+\rightarrow \mathrm{R}_+$ as in Theorem 1 such that (\ref{a8}) holds,
then the trivial solution of (\ref{a1}) is stochastically finite-time stable, and the stochastic settling time $\tau_{x_0}$ satisfies (\ref{a9.1}).

\textbf{Proof:} The proof is similar to that of Theorem 1, and is omitted here.

Thus we may conclude from the theorem the following corollary similar to Corollary 1. 

\textbf{Corollary 2:}  For system (\ref{a1}), suppose that there exists a $C^2$ positive definite and radially unbounded function $U:\mathrm{R}^n\rightarrow\mathrm{R}_+$ satisfying (\ref{a41}). Suppose moreover that there exists a $C^2$ positive definite function $V:\mathrm{R}^n\rightarrow\mathrm{R}_+$, a positive constant $c>0$, such that (\ref{a30}) holds,
then the trivial solution of (\ref{a1}) is stochastically finite-time stable, and the stochastic settling time $\tau_{x_0}$ satisfies (\ref{a30.1}).

\textbf{Remark 9:} The advantage of using multiple Lyapunov functions is that the constraint about $\mathcal{L}V$ can be relaxed without requiring $\mathcal LV\leq 0$ as in Theorem 1, even though in the case of $\mathcal{L}V>0$ (i.e., $\mathcal{L}V$ is a positive definite function) the condition (\ref{a30}) maybe still valid. The example below states this point.

\textbf{Example 3:} Consider the one-dimensional stochastic nonlinear
system,
\begin{eqnarray} \label{a46}
dx=-\frac{1}{2}c_0^2 x^{\frac{p}{q}}dt+c_0x^{\frac{p+q}{2q}}dB(t), \quad x_0\neq0,
\end{eqnarray}
where $p<q$ are positive odd numbers, $c_0\neq0$. Choosing a $C^2$ Lyapunov function $U(x)=x^2$, one knows that $\mathcal{L}U(x)=0$. For any $C^2$ Lyapunov function $V(x)=|x|^{\alpha}$ with $\alpha>2$, it is not hard to verify that
\begin{eqnarray}
\mathcal{L}V(x)=\frac{1}{2}c_0^2(\alpha^2-2\alpha)|x|^{\alpha+\frac{p}{q}-1}>0,
\end{eqnarray}
i.e., $\mathcal{L}V$ is  positive definite. However, for Lyapunov function $V(x)=|x|^{\alpha}$, one can verify that condition (\ref{a30}) is satisfied with $\gamma=1-\frac{q-p}{\alpha q}$ and $c\leq\frac{p+q}{2q}c_0^2\alpha$, clearly $0<\gamma<1$. So, system (\ref{a46}) is stochastically finite-time stable by Corollary 2.


\section {Simulation examples}
\begin{figure}
	\begin{center}
		\epsfig{file=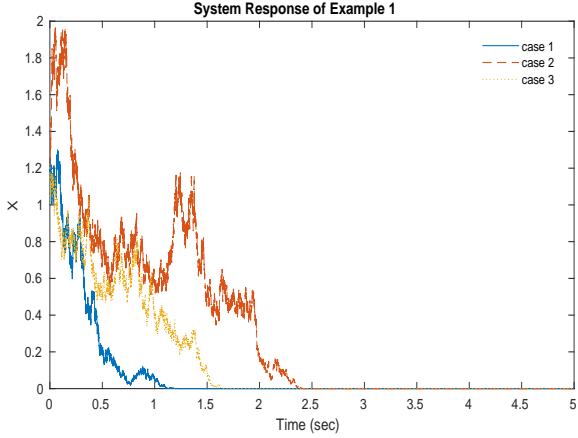,height=6cm, width=9cm}
		\caption{System response of Example 1.}
		\label{fig 1}
	\end{center}
\end{figure}

In this section, we consider Examples 1-3 again, and give their simulation results to illustrate the theory analysis.

 For three cases of Example 1, the initial condition is set to be $x(0)=1.2$. In Case 1, we choose parameters $c_1=0$, $\beta_1=\frac{1}{3}$, $\beta_2=\frac{\beta_1+1}{2}=\frac{2}{3}$, $c_2=\frac{1}{2}$, $c_3=1$ and $\alpha=2$.
	In Case 2, we choose parameters $c_1=1$, $p=1$, $\beta_1=\frac{1}{3}$, $\beta_2=\frac{\beta_1+1}{2}=\frac{2}{3}$, $c_2=\frac{1}{2}$, $c_3=1$ and $\alpha=2$.
In Case 3, we choose parameters $c_1=1$, $p=3$, $\beta_1=\frac{1}{3}$, $\beta_2=\frac{\beta_1+1}{2}=\frac{2}{3}$, $c_2=\frac{1}{2}$, $c_3=1$ and $\alpha=2$. Fig.1 shows the corresponding simulation results. Fig.2 gives the simulation result of Example 2 with initial conditions $x_1(0)=1.5$, $x_2(0)=5$, and Fig.3 shows the simulation result of Example 3 with concrete parameters $c_0=1$, $p=3$ and $q=5$, where the initial condition is set to be $x(0)=2$.

The simulation results clearly show that the trajectories of the corresponding stochastic systems converge rapidly to the equilibrium state in finite time for any given initial values, and verify the effectiveness of theoretical results.

\begin{figure}
	\begin{center}
		\epsfig{file=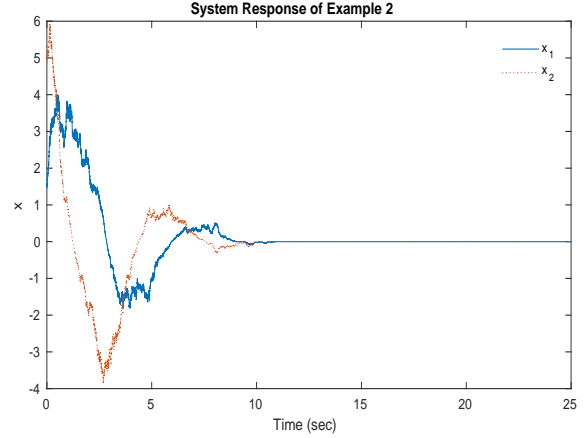,height=6cm, width=9cm}
		\caption{System response of Example 2.}
		\label{fig 1}
	\end{center}
\end{figure}

\begin{figure}
	\begin{center}
		\epsfig{file=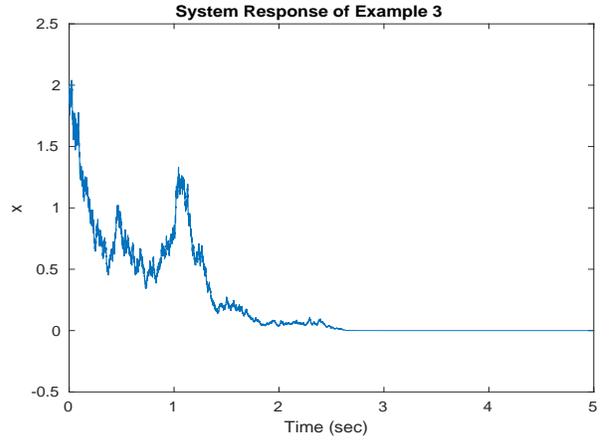,height=6cm, width=9cm}
		\caption{System response of Example 3.}
		\label{fig 1}
	\end{center}
\end{figure}

\section{Conclusion}
In this paper, some new Lyapunov criteria of stochastic finite-time stability are given. Compared with the existing results about stochastic finite-time stability, these new Lyapunov criteria not only relax the constraint on the infinitesimal generator $\mathcal{L}V$,  but also reveal the important role of white-noise in finite-time stabilizing the system. Some examples are constructed to show that these new Lyapunov criteria enable us to construct Lyapunov functions
more easily in applications.


\end{document}